\documentclass[10pt]{amsart}
\usepackage{latexsym, amsmath,amssymb}
 \numberwithin{equation}{section}

\renewcommand{\epsilon}{\varepsilon}
\newtheorem{theorem}{Theorem}

\newtheorem{lemma}[theorem]{Lemma}
\newtheorem{corr}[theorem]{Corollary}

\newtheorem{proposition}[theorem]{Proposition}
\newtheorem{deff}[theorem]{Definition}

\newcommand{\bth}{\begin{theorem}}
\newcommand{\ble}{\begin{lemma}}
\newcommand{\bcor}{\begin{corr}}

\newcommand{\bdeff}{\begin{deff}}

\newcommand{\bprop}{\begin{proposition}}
\newcommand{\ele}{\end{lemma}}
\newcommand{\ecor}{\end{corr}}
\newcommand{\edeff}{\end{deff}}

\numberwithin{theorem}{section}

\newcommand{\eprop}{\end{proposition}}

\renewcommand{\Pi}{\varPi}

\renewcommand{\epsilon}{\varepsilon}

\begin{document}

\title[A pointwise inequality for fractional laplacians]
{A pointwise inequality for fractional laplacians}
\author[A. C\'ordoba]{Antonio C\'ordoba}
\address{Instituto de Ciencias Matem\'aticas CSIC-UAM-UC3M-UCM -- Departamento de Matem\'aticas (Universidad Aut\'onoma de Madrid), 28049 Madrid, Spain} 
\email{antonio.cordoba@uam.es}

\author[A. Mart\'inez]{\'Angel D. Mart\'inez}
\address{Instituto de Ciencias Matem\'aticas CSIC-UAM-UC3M-UCM -- Departamento de Matem\'aticas (Universidad Aut\'onoma de Madrid), 28049 Madrid, Spain} 
\email{angel.martinez@icmat.es}

\begin{abstract}
The fractional laplacian is an operator appearing in several evolution models where diffusion coming from a L\'evy process is present but also in the analysis of fluid interphases. We provide an extension of a pointwise inequality that plays a r\^ole in their study. We begin recalling two scenarios where it has been used. After stating the results, for fractional Laplace-Beltrami and Dirichlet-Neumann operators, we provide a sketch of their proofs, unravelling the underlying principle to such inequalities. 
\end{abstract}
\maketitle



\section{Introduction}

In this exposition we shall extend to a more general framework the following remarkable pointwise inequality
\begin{equation}\label{ineq}
\Lambda^{\alpha}(\phi(f))(x)\leq\phi'(f(x))\cdot\Lambda^{\alpha}f(x)
\end{equation}
valid for any convex function $\phi\in\mathcal{C}^1(\mathbb{R})$ and $f$ in the Schwartz class $\mathcal{S}(\mathbb{R}^n)$, where $\Lambda=(-\Delta)^{\alpha/2}$ is defined as usual through the Fourier transform as
\[\widehat{\Lambda^{\alpha}\theta}(\xi)=|\xi|^{\alpha}\hat{\theta}(\xi)\]
The inequality holds also in the periodic setting and the proof provided in \cite{CC2} follows directly from the representation
\[\Lambda^{\alpha}\theta(x)=c_{n,\alpha}\int_{\mathbb{R}^n}\frac{\theta(x)-\theta(y)}{|x-y|^{n+\alpha}}dy\]
where $c_{n,\alpha}>0$ together with the convexity hypothesis $\phi(f(x))-\phi(f(y))\leq\phi'(f(x))(f(x)-f(y))$. Despite its apparent simplicity its validity is quite surprising given the non-local character of the involved operators. It had also found several interesting applications to non-lineal non-local evolution problems. Let us describe briefly two of them.

\subsection{Transport equations}

Suppose that we are confronted in $\mathbb{R}^n$ or $\mathbb{T}^n$ with the following initial value problem
\[\left\{\begin{array}{l}
\theta_t(x,t)+u(x,t)\cdot\nabla_x\theta(x,t)=-\kappa\Lambda^{\alpha}\theta(x,t)\\
\theta(x,0)=\theta_0(x)
\end{array}\right.\]
where the velocity vector $u$ is divergence-free and $\kappa>0$ is the viscosity coefficient. Then the pointwise inequality \ref{ineq} is crucial to obtain the following maximum principle (see \cite{CC} for details):
\[\|\theta(\cdot,t)\|_p\leq\frac{\|\theta_0\|_p}{(1+C\delta t\|\theta_0\|_p^{p\delta})^{1/p\delta}}\]
where $\delta=\frac{\alpha}{2(p-1)}$, $C=C(\kappa,\alpha,\|\theta_0\|_p)>0$ and $1<p<\infty$.

An important case of this is the surface quasi-geostrophic equation in which case the velocity field of the active scalar $\theta$ is given in terms of the Riesz transforms as $u=(-R_2\theta,R_1\theta)$. In \cite{CC} the pointwise inequality was used to prove that the system above has solution valid in all time $t>0$ for any initial datum $\theta_0$ in the Sobolev space $H^1(\mathbb{R}^2)$. It also appears, for example, in \cite{CV} where regularity results for critical diffusion ($\alpha=1$) are settled. 

\subsection{Interface evolution}

In the theory of water waves or the Muskat and Hele-Shaw flows inside a porous media the inequality \ref{ineq} had played a r\^ole in the study of the evolution of the free boundary between the fluids. In both cases there is a curve (in the plane $\mathbb{R}^2$) or a surface (in the space $\mathbb{R}^3$) whose time evolution has to be controlled. From the fundamental laws (Bernoulli's law for water waves; Darcy's law for porous media) one can asign a velocity field to the moving boundary, obtaining a closed system of rather complicated differential equations.

However, it was early discovered that with such generality the problem is not well-possed. In the case of porous media first Rayleigh, and later Taylor, observed that the linearized problem was unstable if a certain quantity $\sigma$ becomes negative. Therefore, to have a consistent theory one needs to control the evolution of $\sigma$, assuming only its positivity at the initial time.

Choosing and adequate (isothermal) parametrization of the free boundary $\psi:\mathbb{R}^n\rightarrow\mathbb{R}^{n+1}$ (where $n=1,2$) one is lead to estimate the evolution of the Sobolev norms $\frac{d}{dt}\|\psi\|_{H^k}^2$. It turns out that after certain algebraic manipulations the task is finally reduced to the following estimate:
\[\frac{d}{dt}\|\Lambda^k\psi\|_{L^2}^2\leq-\sum_{j=1}^n\sum_{|\alpha|=k}\int\sigma(x,t)D_x^{\alpha}\psi_j(x,t)\Lambda D_x^{\alpha}\psi_j(x,t)dx+O(\|\Lambda^k\psi\|_{L^2}^p)\]
for some positive power $p$. Then \ref{ineq} together with the positivity of the Rayleigh-Taylor term $\sigma$ allows us to dispose of the only dangerous term in the above inequality so as to conclude a well-posedness result. Details might be found in \cite{CCF}.

\section{Statement of results}

The two applications so far considered correspond to models which are both isotropic and homogeneous. However, when those hypothesis are not satisfied it becomes convenient to extend the validity of the pointwise inequality \ref{ineq}. Here we will consider two such extensions, namely to the case of a compact Riemannian manifold $(M,g)$ as well as to the Dirichlet-Neumann operator $\mathcal{D}$ on domains $\Omega\subseteq\mathbb{R}^n$. The conection with the latter class of operator resides in the well-known fact that $\mathcal{D}$ in the upper halfspace in $\mathbb{R}^{n+1}$ coincides precisely with $(-\Delta)^{\frac{1}{2}}$ in $\mathbb{R}^n$ as a simple calculation using Fourier analysis shows.

Let us introduce some notation first, if we denote the metric on $M$ by $|x|_g^2=\sum_{j,k}g_{jk}(x)dx_jdx_k$ recall that the Laplace-Beltrami operator associated with it is given by
\[\Delta_g=\frac{1}{\sqrt{|g|}}\sum_{j,k}\frac{\partial}{\partial x_j}\left(\sqrt{|g|}g^{jk}\frac{\partial}{\partial x_k}\right)\]
where $(g^{jk})=(g_{jk})^{-1}$ and $dx$ denotes the associated volume form as usual. Then the eigenvalues of $-\Delta_g$ are non-negative, numerable and one can find a basis given by the corresponding eigenfunctions $\{\phi_k\}$. The fractional powers $\Lambda_g^{\alpha}$ of $-\Delta_g$, $0\leq\alpha\leq 2$ can be described spectrally as the linear operator that satisfies $\Lambda_g^{\alpha}\phi_k=\lambda_k^{\alpha/2}\phi_k$ for any $k$. Now we can state the first result

\begin{theorem} \label{one}
Given $\alpha\in(0,2]$ and a convex function $\phi\in\mathcal{C}^1(\mathbb{R})$ the following pointwise inequality holds
\[\Lambda_g^{\alpha}(\phi(f))(x)\leq\phi'(f(x))\Lambda^{\alpha}_gf(x)\]
for any $f\in\mathcal{C}^{\infty}(M)$.
\end{theorem}

Given a domain $\Omega\subseteq\mathbb{R}^n$ with smooth ($\mathcal{C}^2$) boundary one may consider the Dirichlet-Neumann operator $\mathcal{D}$ acting on smooth functions $f$ in the boundary $\partial \Omega$. The result reads as follows

\begin{theorem} \label{two}
The following pointwise inequality holds
\[\frac{1}{2m}\mathcal{D}(f^{2m})(x)\leq f(x)^{2m-1}\mathcal{D}f (x)\]
for any positive integer $m\geq 1$.
\end{theorem}

\section{Proof of theorem \ref{one}}

Let us consider the following initial value problems
\[\left\{\begin{array}{l}
\frac{\partial}{\partial t}u(x,t)+\Lambda_g^{\alpha}u(x,t)=0\\
u(x,0)=f(x)
\end{array}\right.\]
and
\[\left\{\begin{array}{l}
\frac{\partial}{\partial t}v(x,t)+\Lambda_g^{\alpha}v(x,t)=0\\
v(x,0)=\phi(f(x))
\end{array}\right.\]
The solutions admit representations
\[u(x,t)=\int_Mf(y)G_{\alpha}(x,y,t)dy\]
and
\[v(x,t)=\int_M\phi(f(y))G_{\alpha}(x,y,t)dy\]
respectively, where the kernel is given by
\begin{equation}\label{kernel}
G_{\alpha}(x,y,t)=\sum_{k}e^{-\lambda_k^{\alpha/2}t}\phi_k(x)\overline{\phi_k(y)}
\end{equation}
Observe that our searched inequality will be an inmediate consequence of the following estimate
\[\frac{\partial}{\partial t}\bigg(v(x,t)-\phi(u(x,t))\bigg)\bigg|_{t=0}\geq 0\]
since $(v-\phi(u))|_{t=0}=0$ it is enough to show that $v(x,t)-\phi(u(x,t))\geq 0$ for any $x\in M$ and $t>0$. But this is a consequence of the positivity nature of the fractional heat kernels $G_{\alpha}(x,y,t)\geq 0$, $\int_MG_{\alpha}(x,y,t)dy=1$ and applying Jensen's inequality together with the convexity hypothesis about $\phi$ one gets the desired inequality
\begin{eqnarray*}
\phi(u(x,t))&=&\phi\left(\int_Mf(y)G_{\alpha}(x,y,t)dy\right)\\
&\leq&\int_M\phi(f(y))G_{\alpha}(x,y,t)dy=v(x,t)
\end{eqnarray*}

To close our argument let us mention that the positivity of $G_{\alpha}$ has been considered among others by S. Bochner (see \cite{B1}; cf. \cite{B2}, \cite{B3}) using the so-called subordination principle. For the sake of completeness we sketch in the following the main lines of its proof.

The first step consists in proving a maximum principle. Suppose that a continuous function $f$ in $M$ is the initial data for the heat equation $\frac{\partial}{\partial t}u+\Delta_gu=0$ then for some positive constant $c_m>0$
\begin{eqnarray*}
\frac{\partial}{\partial t}\|u(\cdot,t)\|^{2m}_{L^{2m}(M)}&=&2m\int_Mu^{2m-1}(x,t)\frac{\partial}{\partial t}u(x,t)dx\\
&=&-c_m\int_Mu^{2m-1}(x,t)\Delta_gu(x,t)dx\\
&=&-c_m\int_Mu^{2m-2}(x,t)|\nabla_gu(x,t)|^2dx\leq 0
\end{eqnarray*}
From which one concludes $\|u(\cdot,t)\|_{L^{2m}}\leq\|f\|_{L^{2m}}$ for any $m\geq 1$, taking limits as $m$ tends to infinity we obtain
\[\|u(\cdot,t)\|_{L^{\infty}(M)}\leq\|f\|_{L^{\infty}(M)}\]
From this estimate we can deduce that $G_2(x,y,t)\geq 0$ by the following argument: assume $G_2(x_0,y_0,t_0)<0$ the, since $\int_MG_2(x,,t)dy=1$ there would be an open set $U\subseteq M$ where
\[\int_UG_2(x_0,y,t_0)dy>1\]
Let $\psi$ be a smooth bump function supported in $U$ such that $0\leq\psi\leq 1$ and such that it approximates the indicator function of $U$ in such a way that
\[\int_M\psi(y)G_2(x_0,y,t_0)dy>1\]
which contradicts the aforementioned maximum principle.

Now we are in position to extend the positivity of the kernel to the rest of values $\alpha\in(0,2)$ for which we will use \ref{kernel} and invoke the Hausdorff-Bernstein-Widder theorem that characterizes representability of functions as a Laplace-Stieltjes transform that assures
\[e^{-\lambda_k^{\alpha}t}=\int_0^{\infty}e^{-\lambda_kt^{1/\alpha}s}d\gamma_{\alpha}(s)\]
where $\gamma_{\alpha}$ is a non-decreasing monotone function. Using this and \ref{kernel} one is tempted to inmediately write
\[G_{\alpha}(x,y,t)=\int_0^{\infty}G_2(x,y,t^{1/{\alpha}}s)d\gamma_{\alpha}(s)\]
from which the positivity of the fractional heat kernel $G_{\alpha}$ would follow. To make this a rigorous argument one needs to control the $L^{\infty}$-norm of the eigenfnctions of $-\Delta_g$ in terms of the corresponding eigenvalues (any polynomial growth would be enough; cf. \cite{SS}) and use Weyl's law. For our purposes the particular case
\[(-\Delta_g)(\phi_k^2)(x)\leq\frac{1}{2}\phi_k(x)\cdot(-\Delta_g\phi_k)=\frac{\lambda_k}{2}\phi_k^2(x)\]
together with the positivity of $(-\Delta_g)^{-1}$, the Sobolev embedding theorem and interpolation provides $\|\phi_k\|_{L^{\infty}}=O(\lambda_k^{n/2})$.

\section{Proof of theorem \ref{two}}

Following the lines of the former proof we propose the following Dirichlet problems in the domain
\[\left\{\begin{array}{lc}
\Delta u=0&\textrm{in $\Omega$}\\
u=f&\textrm{in $\partial \Omega$}
\end{array}\right.\]
and
\[\left\{\begin{array}{lc}
\Delta v=0&\textrm{in $\Omega$}\\
v=f^{2m}&\textrm{in $\partial \Omega$}
\end{array}\right.\]
Then $w=u^{2m}-v$ satisfies
\[\left\{\begin{array}{lc}
\Delta w=2m(2m-1)|\nabla u|^2u^{2m-2}&\textrm{in $\Omega$}\\
w=0&\textrm{in $\partial\Omega$}
\end{array}\right.\]
that is, the subharmonic function $w$ must be non-positive in $\Omega$ since it vanishes at $\partial\Omega$. Consequently, Hopf's lemma implies that $\frac{\partial}{\partial\nu}w(x)>0$ for any $x\in\partial\Omega$ where $\nu$ is the exterior normal to the domain $\Omega$. But this is exactly the desired inequality.

\section{Acknowledgments}

The first author is indebted to E. M. Stein for asking extensions of the inequality \ref{ineq} to the compact riemannian manifold setting, suggesting an alternative approach. Partially supported by MTM2011-2281 project of the MCINN (Spain).


\begin{thebibliography}{10}
\bibitem{B1} 
Bochner, S., {\em Quasi-Analytic Functions, Laplace Operator, Positive Kernels}, Ann. of Math. Vol. 51, No. 1 (1950), pp. 68-91.
\bibitem{B2} 
Bochner, S., {\em Diffusion Equation and Stochastic Processes}, Proc. Nat. Aca. Sci. USA, Vol. 35. No. 7 (1949), pp. 368-370.
\bibitem{B3} 
Bochner, S., {\em Stable laws of probability and completely monotone functions}, Duke Math. J. Vol. 3. No. 4 (1937), pp. 726-728.
\bibitem{CC} 
C\'ordoba, A.; C\'ordoba, D., {\em A Maximum Principle Applied to Quasi-Geostrophic Equations}, Commun. Math. Phys. 249 (2004), pp. 511-528.
\bibitem{CC2} 
C\'ordoba, A.; C\'ordoba, D., {\em A pointwise estimate for fractionary derivatives with applications to partial differential equations.} Proceedings of the National Academy of Sciences of the United States of America 100 (26) (2003), pp. 15316-15317.
\bibitem{CCF} 
C\'ordoba, A.; C\'ordoba, D.; Gancedo, F., {\em Interface evolution: the Hele-Shaw and Muskat problems}, Ann. of Math. Vol. 173 (2011), pp. 477-542.
\bibitem{CV} 
Caffarelli, L. A.; Vasseur, A., {\em Drift diffusion equations with fractional diffusion and the quasi-geostrophic equation}, Ann. of Math. Vol. 171, No. 3 (2010), pp. 1903-1930.
\bibitem{SS} 
Seeger, A.; Sogge, C. D., {\em Bounds for eigenfunctions of differential operators}, Indiana U. Math. J., Vol. 38 (1989), No. 38.
\end{thebibliography}
\end{document}